\documentclass[12pt]{article}
\usepackage{setspace}
\RequirePackage[OT1]{fontenc}
\usepackage{booktabs}
\usepackage{amsmath,amsfonts,amssymb,amsthm}
\usepackage{graphicx}

\usepackage{paralist}
\usepackage{amsbsy}

\usepackage[mathscr]{eucal}

\usepackage{bm}

\usepackage{xspace}

\mathchardef\mhyphen="2D

\newcommand{\Real}{\mathbb{R}}

\newcommand{\V}[1]{{\bm{\mathbf{\MakeLowercase{#1}}}}} 
\newcommand{\VE}[2]{\MakeLowercase{#1}_{#2}} 

\newtheorem{theorem}{Theorem}[section]
\newtheorem{proposition}[theorem]{Proposition}

\begin{document}
\title{A Look at the Generalized Heron Problem through the Lens of Majorization-Minimization}
\author{Eric C. Chi and Kenneth Lange}
\date{}
\maketitle

\begin{abstract}
In a recent issue of this journal, Mordukhovich et al.\ pose and solve an interesting non-differentiable generalization of the Heron problem in the framework of modern convex analysis. In the generalized Heron problem one is given $k+1$ closed convex sets in $\Real^d$ equipped with its Euclidean norm and asked to find the point in the last set such that the sum of the distances to the first $k$ sets is minimal.  In later work the authors generalize the Heron problem even further, relax its convexity assumptions, study its theoretical properties, and pursue subgradient algorithms for solving the convex case.  Here, we revisit the original problem solely from the numerical perspective.  By exploiting the majorization-minimization (MM) principle of computational statistics and rudimentary techniques from differential calculus,  we are able to construct a very fast algorithm for solving the Euclidean version of the generalized Heron problem.
\end{abstract}

\section{Introduction.}

In a recent article in this journal, Mordukhovich et al.~\cite{Mordukhovich2012} presented the following generalization of the classical Heron problem. Given a collection of closed convex sets $\{C_1, \ldots, C_k\}$ in $\Real^d$, find a point \V{x} in the closed convex set $S \subset \Real^d$ such that the sum of the Euclidean distances from $\V{x}$ to $C_1$ through $C_k$ is minimal. In other words,

\begin{equation}
\label{eq:Heron}
\text{minimize $D(\V{x}) := \sum_{i=1}^k d(\V{x}, C_i)$ subject to $\V{x} \in S$},
\end{equation}
where $d(\V{x}, \Omega) = \inf\{ \lVert \V{x} - \V{y} \rVert : \V{y} \in \Omega \}$.

A rich history of special cases motivates this problem formulation. When $k = 2$, $C_1$ and $C_2$ are singletons, and $S$ is a line, we recover the problem originally posed by the ancient mathematician Heron of Alexandria. The special case where $k=3$; $C_1$, $C_2,$ and $C_3$ are singletons; and $S = \Real^2$ was suggested by Fermat nearly 400 years ago and solved by Torricelli \cite{Krarup1997}.  In his {\em Doctrine and Application of Fluxions}, Simpson generalized the distances to weighted distances. In the 19th century, Steiner made several fundamental contributions, and his name is sometimes attached to the problem \cite{Courant1961,Gueron2002}.  At the turn of the 20th century, the German economist Weber generalized Fermat's problem to an arbitrary number of singleton sets $C_i$.  Weiszfeld published the first iterative algorithm\footnote{Kuhn \cite{Kuhn1973} points out that Weiszfeld's algorithm has been rediscovered several times.} for solving the Fermat-Weber problem in 1937 \cite{Weiszfeld1937, Weiszfeld2009}. In the modern era, the Fermat-Weber problem has enjoyed a renaissance in various computational guises. Both the problem and associated algorithms serve as the starting point for many advanced models in location theory \cite{Love1988, Wesolowsky1993}.

The connections between celebrated problems such as the Fermat-Weber problem and the generalized Heron problem were noted earlier by Mordukhovich et al.~\cite{Mordukhovich2011b}.  In subsequent papers \cite{Mordukhovich2011,Mordukhovich2011b}, they generalize the Heron problem further to arbitrary closed sets, $C_1, \ldots, C_k$ and $S$ in a Banach space. Readers are referred to their papers for a clear treatment of how one solves these abstract versions of the generalized Heron problem with state-of-the-art tools from variational analysis. 

Here we restrict our attention to the special case of Euclidean distances presented by Mordukhovich et al.~\cite{Mordukhovich2011b}.  Our purpose is take a second look at this simple yet most pertinent version of the  problem from the perspective of algorithm design. Mordukhovich et al. \cite{Mordukhovich2011, Mordukhovich2012, Mordukhovich2011b} present an iterative subgradient algorithm for numerically solving problem (\ref{eq:Heron}) and its generalizations, a robust choice when one desires to assume nothing beyond the convexity of the objective function. Indeed, the subgradient algorithm works if the Euclidean norm is exchanged for an arbitrary norm. However, it is natural to wonder if there might be better alternatives for the finite-dimensional version of the problem with Euclidean distances. Here we present one that generalizes Weiszfeld's algorithm by invoking the majorizati\-on-minimization (MM) principle from computational statistics.  Although the new algorithm displays the same kind of singularities that plagued Weiszfeld's algorithm \cite{Kuhn1973}, the dilemmas can be resolved by slightly perturbing problem (\ref{eq:Heron}), which we refer to as the generalized Heron problem for the remainder of this article.  In the limit, one recovers the solution to the unperturbed problem.  As might be expected, it pays to exploit special structure in a problem. The new MM algorithm is vastly superior to the subgradient algorithms in computational speed for Euclidean distances.

Solving a perturbed version of the problem by the MM principle yields extra dividends as well.  The convergence of MM algorithms on smooth problems is well understood theoretically. This fact enables us to show that solutions to the original problem can be characterized without appealing to the full machinery of convex analysis dealing with non-differentiable functions and their subgradients.  Although this body of mathematical knowledge is definitely worth learning, it is remarkable how much progress can be made with simple tools. The good news is that we demonstrate that crafting an iterative numerical solver for problem (\ref{eq:Heron}) is well within the scope of classical differential calculus.  Our resolution can be understood by undergraduate mathematics majors. 

As a brief summary of things to come, we begin by recalling background material on the MM principle and convex analysis of differentiable functions. This is followed with a derivation of the MM algorithm for problem (\ref{eq:Heron}) and consideration of a few relevant numerical examples. We end by proving convergence of the algorithm and characterizing solution points.

\section{The MM Principle.}

Although first articulated by the numerical analysts Ortega and Rheinboldt \cite{Ortega1970}, the MM principle currently enjoys its greatest vogue in computational statistics \cite{Becker1997, Lange2000}. The basic idea is to convert a hard optimization problem (for example, non-differentiable) into a sequence of simpler ones (for example, smooth). The MM principle requires majorizing the objective function $f(\V{y})$ by a surrogate function $g(\V{y} \mid \V{x})$ anchored at the current point $\V{x}$.  Majorization is a combination of the tangency condition $g(\V{x} \mid \V{x}) =  f(\V{x})$ and the domination condition $g(\V{y} \mid \V{x})  \geq f(\V{y})$ for all $\V{y} \in \Real^d$.  The associated MM algorithm is defined by the iterates
\begin{equation}
  \label{eq:generic-MM-iterate}
  \V{x}_{k+1} := \arg \min_{\V{y} \in S} g(\V{y} \mid \V{x}_{k}).
\end{equation}
Because 
\begin{equation}
  f(\V{x}_{k+1}) \leq g(\V{x}_{k+1} \mid \V{x}_{k}) \leq g(\V{x}_{k} \mid \V{x}_{k}) = f(\V{x}_{k}),
\end{equation}
the MM iterates generate a descent algorithm driving the objective function downhill.  Constraint satisfaction is enforced in finding $\V{x}_{k+1}$. Under appropriate regularity conditions, an MM algorithm is guaranteed to converge to a local minimum of the original problem \cite{Lange2010}. 

\section{Background on Convex Analysis.}

As a prelude to deriving an MM algorithm, we review some basic facts from convex analysis in the limited context of differentiable functions. Deeper treatments can be found in the references \cite{Bertsekas2009, Borwein2000, Hiriart-Urruty2004,Rockafellar1996, Ruszczynski2006}. Recall that a differentiable function $f(\V{Y})$ is convex if and only if its domain $S$ is convex and
\begin{equation}
\label{eq:support}
f(\V{y}) \geq f(\V{x}) + \langle \nabla f(\V{x}), \V{y} - \V{x} \rangle,
\end{equation}
for all $\V{x}, \V{y} \in S$. Provided $f(\V{x})$ is twice differentiable, it is convex when its second differential $d^2 f(\V{x})$ is positive semidefinite for all $\V{x}$ and strictly convex when $d^2 f(\V{x})$ is positive definite for all $\V{x}$. These characterizations are a direct consequence of executing a second-order Taylor expansion of $f(\V{y})$ and applying the supporting hyperplane inequality (\ref{eq:support}). The supporting hyperplane inequality (\ref{eq:support}) also leads to a succinct necessary and sufficient condition for a global minimum.  A point $\V{x} \in S$ is a global minimizer of $f(\V{y})$ on $S$ if and only if 
\begin{equation}
\langle \nabla f(\V{x}), \V{y} - \V{x} \rangle \geq 0
\end{equation}
for all $\V{y} \in S$. Intuitively speaking, every direction pointing into $S$ must lead uphill.

We conclude this section by reviewing projection operators \cite{Lange2010}. Denote the projection of $\V{x}$ onto a set $\Omega \subset \Real^d$ by $P_\Omega(\V{x})$. By definition $P_\Omega(\V{x})$ satisfies
\begin{equation*}
P_\Omega(\V{x}) := \underset{\V{y} \in \Omega}{\arg\min} \lVert \V{x} - \V{y} \rVert.
\end{equation*}
If $\Omega$ is a closed convex set in $\Real^d$, then $P_\Omega(\V{x})$ exists and is unique. Furthermore, the projection operator is non-expansive in the sense that
\begin{equation*}
\lVert P_\Omega(\V{x}) - P_\Omega(\V{y}) \rVert \leq \lVert \V{x} - \V{y} \rVert
\end{equation*}
for all $\V{x}, \V{y} \in \Real^d$. Non-expansion clearly entails continuity.   Explicit formulas for the projection operator $P_\Omega(\V{x})$ exist when $\Omega$ is a box, Euclidean ball, hyperplane, or halfspace. Fast algorithms for computing $P_\Omega(\V{x})$ exist for the unit simplex, the $\ell_1$ ball, and the cone of positive semidefinite matrices \cite{Duchi2008, Michelot1986}. 

The projection operator and the distance function are intimately related through the gradient identity $\nabla d(\V{x}, C)^2 = 2 [\V{x} - P_C(\V{x})]$. A standard proof of this fact can be found in reference \cite[p.~181]{Hiriart-Urruty2004}. If $d(\V{x}, C)^2>0$, then the chain rule gives
\begin{equation*}
\nabla d(\V{x}, C)= \nabla \sqrt{d(\V{x},C)^2} = \frac{\V{x} - P_C(\V{x})}{d(\V{x}, C)}.
\end{equation*}
On the interior of $C$, it is obvious that $\nabla d(\V{x}, C)= {\bf 0}$. In contrast, differentiability of $d(\V{x}, C)$ at boundary points of $C$ is not guaranteed.

\section{An MM Algorithm for the Heron Problem.}

Since it adds little additional overhead, we recast problem (\ref{eq:Heron}) in the Simpson form
\begin{equation}
\label{eq:ourHeron}
\text{minimize $D(\V{x}) := \sum_{i=1}^k \gamma_i d(\V{x}, C_i)$ subject to $\V{x} \in S$}
\end{equation}
involving a convex combination of the distances $d(\V{x}, C_i)$ with positive weights $\gamma_i$ as suggested in \cite{Mordukhovich2011b}. We first derive an MM algorithm for solving problem (\ref{eq:ourHeron}) when $S \cap C_i = \emptyset$ for all $i$. This exercise will set the stage for attacking the more general case where $S$ intersects one or more of the $C_i$.  In practice quadratic majorization is desirable because it promotes exact solution of the minimization step of the MM algorithm.  It takes two successive majorizations to achieve quadratic majorization in our setting. The first is the simple majorization
\begin{equation*}
d(\V{x}, C_i) \leq \lVert \V{x} - P_{C_i}(\V{x}_m) \rVert 
\end{equation*}
flowing directly from the definition of the distance function. The second is the majorization
\begin{equation}
\label{eq:MM_univariate}
\sqrt{u} \leq \sqrt{u_m} + \frac{1}{2\sqrt{u_m}} (u - u_m),
\end{equation}
of the concave function $\sqrt{u}$ on the interval $(0,\infty)$. The combination of these two majorizations yields the quadratic majorization
\begin{equation}
\label{eq:quad_majorization}
d(\V{x}, C_i) \le 
\lVert \V{x}_m - P_{C_i}(\V{x}_m) \rVert +
\frac{\lVert \V{x} - P_{C_i}(\V{x}_m) \rVert^2 - \lVert \V{x}_m - P_{C_i}(\V{x}_m) \rVert^2}
{2 \lVert \V{x}_m - P_{C_i}(\V{x}_m) \rVert}.
\end{equation}
Summing these majorizations over $i$ leads to quadratic majorization of $D(\V{x})$ and ultimately to the MM algorithm map
\begin{equation*}
\psi(\V{x}) = \underset{\V{z} \in S}{\arg\min} \left\{ \frac{1}{2} \sum_{i=1}^k \VE{w}{i}\lVert \V{z} - P_{C_i}(\V{x}) \rVert^2 
\right \}
\end{equation*}
with weights $\VE{w}{i} = \VE{\gamma}{i}\lVert \V{x} - P_{C_i}(\V{x}) \rVert^{-1}$.  When the $C_i$ are singletons and $S = \Real^d$, the map
$\psi(\V{x})$ implements Weiszfeld's algorithm for solving the Fermat-Weber problem \cite{Weiszfeld1937, Weiszfeld2009}. 

The quadratic majorization of $D(\V{x})$  just derived can be rewritten as
\begin{equation*}
g(\V{x} \mid \V{x}_m) =
\frac{1}{2}\left(\sum_{i=1}^k \VE{w}{i} \right)\Big\| \V{x} - \sum_i \VE{\alpha}{i} P_{C_i}(\V{x}_m) \Big\|^2 + c,
\end{equation*}
where
\begin{equation*}
\alpha_i = \frac{w_i}{\sum_{i=1}^k w_i},
\end{equation*}
and $c$ is a constant that does not depend on $\V{x}$. Thus, the MM update boils down to projection onto $S$ of a convex combination of the projections of the previous iterate onto the sets $C_i$; in symbols
\begin{equation}
\label{eq:MM_update}
\V{X}_{m+1} = P_S\Big[\sum_i \VE{\alpha}{i}  P_{C_i}(\V{x}_m) \Big].
\end{equation}

The majorization (\ref{eq:quad_majorization}) involves dividing by 0 when $\V{x}_m$ belongs to $C_i$. This singularity also bedevils Weiszfeld's algorithm.  Fortunately, perturbation of the objective function salvages the situation. One simply replaces the function $D(\V{x})$ by the related function
\begin{eqnarray*}
D_{\epsilon}(\V{x}) & = & \sum_{j=1}^k \gamma_j \sqrt{d(\V{x},C_j)^2+\epsilon}
\end{eqnarray*}
for $\epsilon$ small and positive. Ben-Tal and Teboulle \cite{Ben-Tal1989} cover further examples of this perturbation strategy. In any case observe that the smooth function $f_{\epsilon}(u)=\sqrt{u^2+\epsilon}$ has derivatives
\begin{eqnarray*}
f_{\epsilon}'(u) & = & \frac{ u }{ \sqrt{u^2+\epsilon}}, \quad 
f_{\epsilon}''(u) \:\;\, = \:\;\, \frac{\epsilon }{ (u^2+\epsilon)^{3/2}}
\end{eqnarray*}
and is therefore strictly increasing and strictly convex on the
interval $[0,\infty)$.  Hence, the function $D_{\epsilon}(\V{x})$ is also convex.
Because $\sqrt{u^2+\epsilon}-\sqrt{\epsilon}$ is a good approximation to $u \ge 0$, the solutions of the two problems should be close. In fact, we will show later that the minimum point of $D_{\epsilon}(\V{x})$ tends to the minimum point of $D(\V{x})$ as $\epsilon$ tends to 0. In the presence of multiple minima, this claim must be rephrased in terms of cluster points.

The majorization $d(\V{x},C_j) \le \|\V{x}-P_{C_j}(\V{x}_m)\|$ around the current iterate $\V{x}_m$ yields the majorization
\begin{eqnarray*}
\sqrt{d(\V{x},C_j)^2+\epsilon} & \le & \sqrt{\|\V{x}-P_{C_j}(\V{x}_m)\|^2+\epsilon} .
\end{eqnarray*}
Application of the majorization (\ref{eq:MM_univariate}) implies the further majorization
\begin{eqnarray*}
D_{\epsilon}(\V{x}) & \le & \frac {1 }{ 2} \sum_{j=1}^k \gamma_j 
\frac{\|\V{x}-P_{C_j}(\V{x}_m)\|^2 }{ \sqrt{\|\V{x}_m-P_{C_j}(\V{x}_m)\|^2+\epsilon}}+c,
\end{eqnarray*}
where $c$ is an irrelevant constant.  The corresponding MM update $\V{x}_{m+1}$ is identical to the previous MM update (\ref{eq:MM_update}) except for one difference. The weights $\VE{w}{i}$ are now defined by the benign formula
\begin{equation*}
\VE{w}{i} = \frac{\gamma_i}{\sqrt{\lVert \V{x}_m - P_{C_i}(\V{x}_m) \rVert^2 + \epsilon}}
\end{equation*}
involving no singularity.

\section{Examples.}

We now consider four examples illustrating the performance of the MM algorithm and framing our expectations for convergence. The subgradient algorithm \cite{Mordukhovich2012} serves as a benchmark for comparison throughout. This algorithm relies on the updates
\begin{equation*}
\V{x}_{m+1} = P_S\Big[\V{x}_m - \eta_m \sum_{i=1}^k \gamma_i \V{v}_{im} \Big],
\end{equation*}
where
\begin{equation*}
\V{v}_{im} = \begin{cases}
\frac{\V{x}_m - P_{C_i}(\V{x}_m)}{d(\V{x}_m, C_i)} & \text{if $\V{x}_m \not\in C_i$} \\
0 & \text{if $\V{x}_m \in C_i$},
\end{cases}
\end{equation*}
and the nonnegative constants $\eta_m$ satisfy $\sum_{m=1}^\infty \eta_m = \infty$ and $\sum_{m=1}^\infty \eta_m^2 < \infty$.  The weights $\gamma_i$ equal $1$ in all examples except the last.

\begin{table}[ht]
\begin{center}
\begin{tabular}{rccc}
  \hline
 Iteration & $x_1$ & $x_2$ & $x_3$ \\ 
  \hline
             1 & 0.00000000000000 & 2.00000000000000 & 0.00000000000000 \\ 
               2 & -0.93546738305698 & 1.66164748416805 & 0.10207032020482 \\ 
              3 & -0.92881282698649 & 1.63915389878166 & 0.08424264751830 \\ 
            4 & -0.92645373003448 & 1.63220797263449 & 0.08007815377225 \\ 
              5 & -0.92567602259658 & 1.63004821970935 & 0.07911751670489 \\ 
              6 & -0.92542515217106 & 1.62937435413374 & 0.07889815178685 \\ 
               7 & -0.92534495711879 & 1.62916364685109 & 0.07884864943702 \\ 
               8 & -0.92531944712805 & 1.62909766226627 & 0.07883765997470 \\ 
          9 & -0.92531135783449 & 1.62907697582185 & 0.07883527888603 \\ 
           10 & -0.92530879826106 & 1.62907048520349 & 0.07883478238381 \\ 
             20 & -0.92530761702316 & 1.62906751412014 & 0.07883466748783 \\ 
             30 & -0.92530761701184 & 1.62906751409212 & 0.07883466748878 \\ 
             50 & -0.92530761701184 & 1.62906751409212 & 0.07883466748878 \\ 
   \hline
\end{tabular}
\end{center}
\caption{Cubes and ball example in $\Real^3$: MM Algorithm.}
\label{tab:MM_R3}
\end{table}

\subsection{Five Cubes and a Ball in $\Real^3$.}

Our first example is taken from the reference \cite{Mordukhovich2012}. This three-dimensional example involves five cubes $C_i$ with side lengths equal to 2 and centers $(0,-4,0), (-4,2,-3), (-3,-4,2)$, $(-5,4,4)$, and $(-1,8,1)$. The set $S$ is a ball with center $(0,2,0)$ and radius 1.  Iteration commences at the point $\V{x}_1 = (0, 2, 0) \in S$ and takes subgradient steps with $\eta_m = 1/m$. Table~\ref{tab:MM_R3} shows the MM iterates with $\epsilon = 0$. Convergence to machine precision occurs within 30 iterations. In contrast Table~\ref{tab:pg_R3} shows that parameter values $(x_1,x_2,x_3)$ are still changing after $10^6$ subgradient iterates.  For brevity we omit a second example of four squares and a disk in $\Real^2$ from the same source \cite{Mordukhovich2012}. In this example the superiority of the MM algorithm over the subgradient algorithm is equally evident.

\begin{table}[ht]
\begin{center}
\begin{tabular}{rccc}
  \hline
  Iteration & $x_1$ & $x_2$ & $x_3$ \\ 
  \hline
              1 & 0.00000000000000 & 2.00000000000000 & 0.00000000000000 \\ 
             10 & -0.92583298353433 & 1.63051788239768 & 0.07947484741743 \\ 
            100 & -0.92531325048300 & 1.62908232435160 & 0.07883822912883 \\ 
           1000 & -0.92530767419684 & 1.62906766065418 & 0.07883468589312 \\ 
          10000 & -0.92530761758555 & 1.62906751554109 & 0.07883466757273 \\ 
         100000 & -0.92530761701755 & 1.62906751410641 & 0.07883466748904 \\ 
       1000000 & -0.92530761701233 & 1.62906751409334 & 0.07883466748881 \\ 
          1500000 & -0.92530761701231 & 1.62906751409328 & 0.07883466748881 \\ 
         2000000 & -0.92530761701229 & 1.62906751409324 & 0.07883466748881 \\ 
   \hline
\end{tabular}
\end{center}
\caption{Cubes and ball example in $\Real^3$:  Subgradient Algorithm.}
\label{tab:pg_R3}
\end{table}

\subsection{The Closest Point to Three Disks in $\Real^2.$}

\begin{figure}
\centering
\includegraphics[scale=0.45]{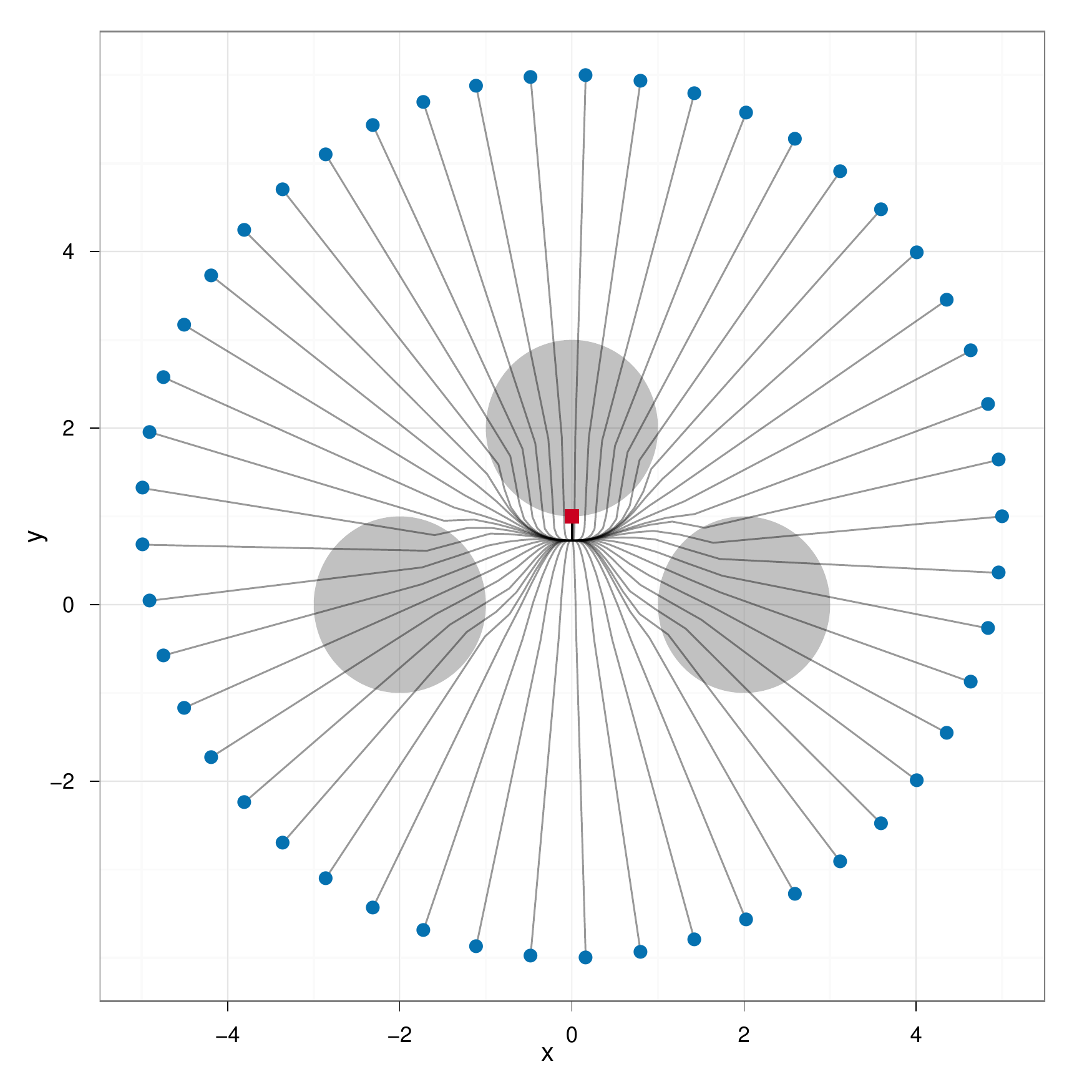}
\caption{Finding the closest point to three disks in $\Real^2.$}
\label{fig:three_disks}
\end{figure}

This example from the reference \cite{Mordukhovich2011} illustrates the advantage of minimizing a sequence of approximating functions $D_{\epsilon_m}(\V{x})$. The sets $C_i$ are three unit balls in $\Real^2$ centered at $(0,2), (2,0)$, and $(-2,0)$. The set $S$ equals $\Real^2$. The minimum distance occurs at $(0,1)$ as can be easily verified by checking the optimality conditions spelled out in Proposition 4.3 in \cite{Mordukhovich2011}. Figure~\ref{fig:three_disks} displays the iteration paths for 50 different starting values (dots) and their corresponding fixed point (the square). Along the $m$th leg of the path we set $\epsilon_m$ to be $\max\{10^{-m},10^{-16}\}$. The solution to the current problem is taken as the initial point for the next problem. All solution paths initially converge to a point just below (0,1) and then march collectively upwards to (0,1). The passage of the MM iterates through the unit balls is facilitated by our strategy of systematically reducing $\epsilon$. Table~\ref{tab:three_disks} shows the subgradient and MM iterates starting from the point (5,7).

\begin{table}
\centering
\begin{tabular}{rccrcc}
\toprule
\multicolumn{3}{c}{Subgradient Algorithm} & \multicolumn{3}{c}{MM Algorithm} \\
\cmidrule(lr){1-3} \cmidrule(lr){4-6}
Iteration & $\VE{x}{1}$ & $\VE{x}{2}$ & Iteration  & $\VE{x}{1}$ & $\VE{x}{2}$ \\
\midrule
 10 & 0.7092649 & 1.2369866  &  10 & 0.2674080 & 0.7570688 \\
 100 & 0.0558764 & 0.9973310  & 100 & 0.0000000 & 0.7249706 \\
 1,000 & 0.0046862 & 0.9993844 & 1,000 &0.0000000 & 0.9998002 \\
 10,000 & 0.0003955 & 0.9999274 & 1,800 & 0.0000000 & 0.9999999 \\
 100,000 & 0.0000334 & 0.9999957 & 1,850 & 0.0000000 & 1.0000000 \\
 1,000,000 & 0.0000028 & 0.9999998 & 1,900 & 0.0000000 & 1.0000000 \\
 \bottomrule
\end{tabular}
\caption{Three disks example in $\Real^2$ starting from (5,7).}
\label{tab:three_disks}
\end{table}

\subsection{Three Collinear Disks in $\Real^2$.}

\begin{figure}
\centering
\includegraphics[scale=0.375]{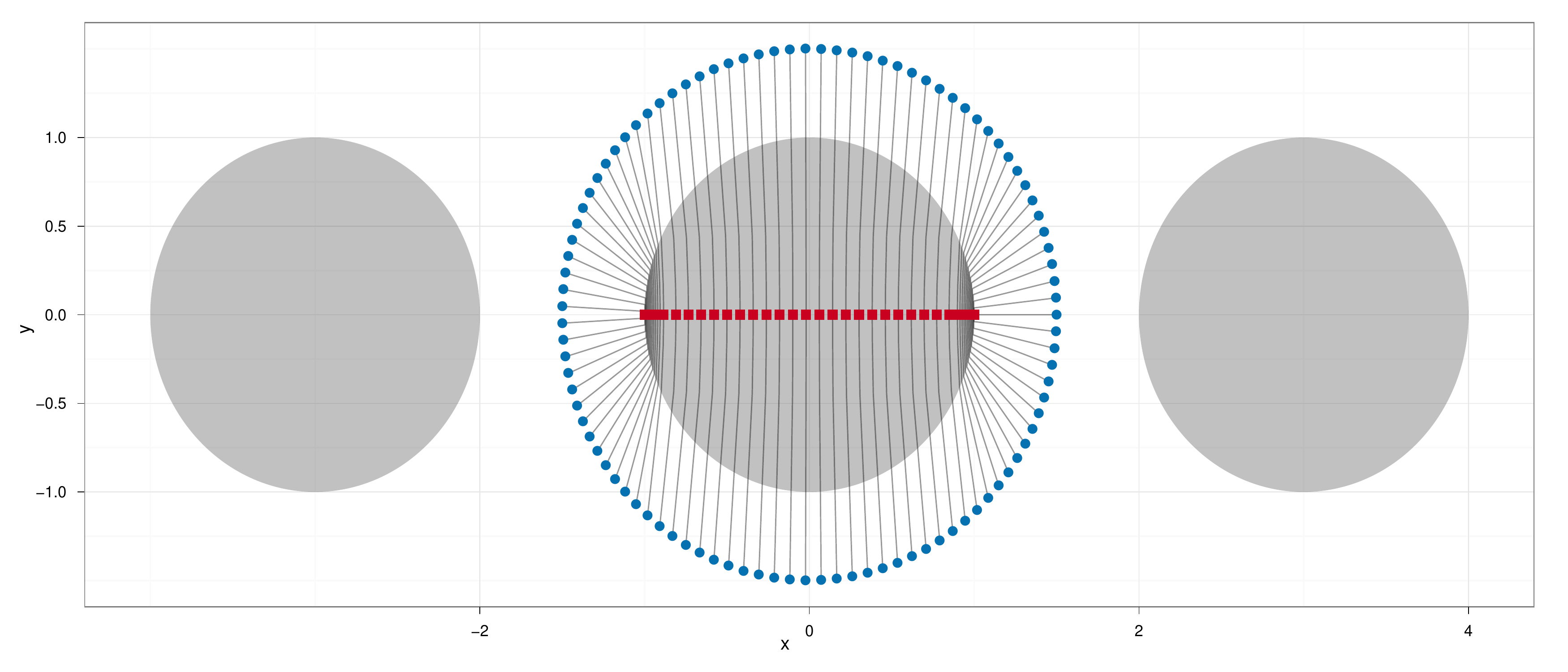}
\caption{An example with a continuum of solutions.}
\label{fig:continuum}
\end{figure}

Here we illustrate the behavior of the MM algorithm when there is more than one solution. Consider two unit balls in $\Real^2$ centered at $(2,0)$, and $(-2,0)$, and take $S$ to be the unit ball centered at the origin. There is a continuum of solutions extending along the line segment from $(-1,0)$ to $(1,0)$, as can be verified by the optimality conditions provided by Theorem 3.2 in \cite{Mordukhovich2012}. Figure~\ref{fig:continuum} shows the iteration paths for 100 different initial values (dots) and their corresponding fixed points (squares). In this example we take $\epsilon = 0$. Although the iterates are not guaranteed to converge and may in principle cycle among multiple cluster points, this behavior is not observed in practice. The iterates simply converge to different fixed points depending on where they start. Table~\ref{tab:three_collinear_disks} compares the iterations for the subgradient method and the MM algorithm starting from the point (1.5,0.25). The two algorithms converge to different solution points but at drastically different rates.

\begin{table}
\centering
\begin{tabular}{rccrcc}
\toprule
\multicolumn{3}{c}{Subgradient Algorithm} & \multicolumn{3}{c}{MM Algorithm} \\
\cmidrule(lr){1-3} \cmidrule(lr){4-6}
Iteration & $\VE{x}{1}$ & $\VE{x}{2}$ & Iteration  & $\VE{x}{1}$ & $\VE{x}{2}$ \\
\midrule
 10,000       & 0.9997648 & 0.0000223  &  10 & 0.9941149 & 0.0001308 \\
 100,000     & 0.9997648 & 0.0000040   & 20 & 0.9941149 & 0.0000000 \\
 1,000,000 & 0.9997648 & 0.0000007 &  30 & 0.9941149 & 0.0000000 \\
 \bottomrule
\end{tabular}
\caption{Three collinear disks example in $\Real^2$ starting from (1.5, 0.25)}
\label{tab:three_collinear_disks}
\end{table}

\subsection{Kuhn's Problem.}
\begin{figure}
\centering
\includegraphics[scale=0.35]{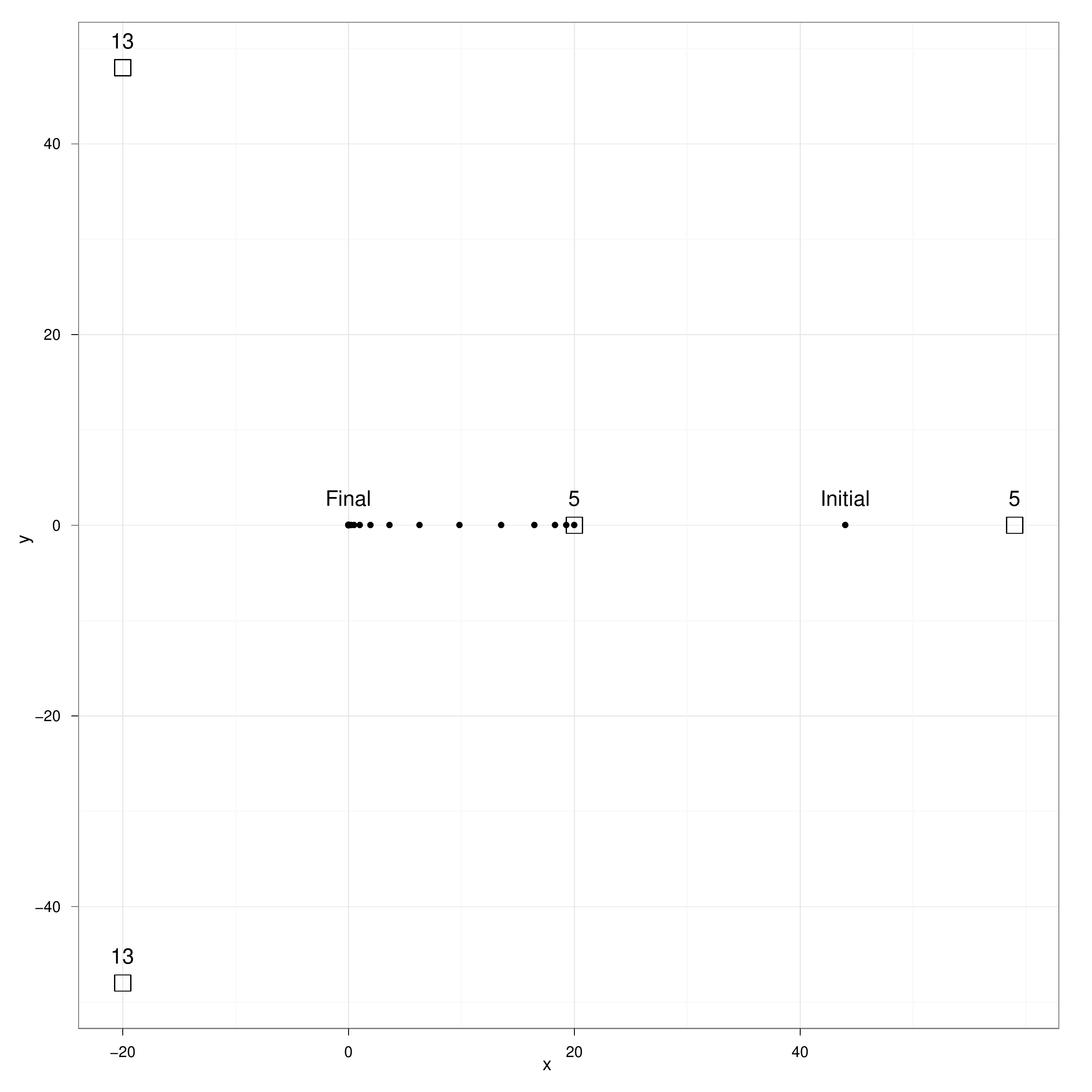}
\caption{A problem where Weiszfeld's algorithm fails to converge.}
\label{fig:Kuhn}
\end{figure}

Our last example was originally concocted by Kuhn \cite{Kuhn1967} to illustrate how Weiszfeld's algorithm can stall when its iterates enter one of the sets $C_i$.  Although this event rarely occurs in practice, characterizing the initial conditions under which it happens has been a subject of intense scrutiny \cite{Brimberg1995, Brimberg2003,Canovas2002,Chandrasekaran1989,Kuhn1973}.  The occasional failure of Weiszfeld's algorithm  prompted Vardi and Zhang  \cite{Vardi2001} to redesign it.  Their version preserves the descent property but differs substantially from ours.  In any event the example shown in Figure~\ref{fig:Kuhn} involves two points with weights $\gamma_i$ proportional to 5 placed at (59,0) and (20,0) and two more points with weights proportional to 13 placed at (-20, 48) and (-20, -48). The optimal point is the origin. Starting at (44,0), Weiszfeld's algorithm stalls at (20,0) after one iteration. Our MM iterates (dots) with $\epsilon$ decreasing from $0.1$ to $0$, in contrast, move across (20.0) and correctly converge to (0,0) to within machine precision in 99 steps. Table~\ref{tab:Kuhn} compares the progress achieved by the MM and subgradient methods. Note that when $\epsilon$ is $0.1$, the MM algorithm overshoots the true answer and then comes back to $(0,0)$ after setting $\epsilon$ to be 0. The subgradient algorithm makes solid progress early but subsequently slows down on this almost smooth problem.

\begin{table}
\centering
\begin{tabular}{rccrcc}
\toprule
\multicolumn{3}{c}{Subgradient Algorithm} & \multicolumn{3}{c}{MM Algorithm} \\
\cmidrule(lr){1-3} \cmidrule(lr){4-6}
Iteration & $\VE{x}{1}$ & $\VE{x}{2}$ & Iteration  & $\VE{x}{1}$ & $\VE{x}{2}$ \\
\midrule
 10      & 8.6984831 & 0.0000000  &  10 & 1.9448925 & 0.0000000 \\
 1,000    & 1.2966354 & 0.0000000   & 30 &  -0.0011998 & 0.0000000 \\
 100,000 & 0.1845171 & 0.0000000 &  60 &  -0.0012011 & 0.0000000 \\
  10,000,000 & 0.0259854 & 0.0000000 & 90 & 0.0000000 & 0.0000000\\ 
 \bottomrule
\end{tabular}
\caption{Kuhn's problem}
\label{tab:Kuhn}
\end{table}

\section{Convergence Theory.}

Before embarking on a proof of convergence, it is prudent to discuss whether a minimum point exists and is unique.  Recall that a continuous function attains its minimum on a compact set.  Thus, problem (\ref{eq:ourHeron}) possesses a minimum whenever $S$ is bounded.  If $S$ is unbounded, then one can substitute boundedness of one or more of the sets $C_i$.  In this circumstance  $D(\V{x})$ is coercive in the sense that $\lim_{\|\V{x}\| \to \infty} D(\V{x}) = \infty$.  As pointed out in Proposition 3.1 of the reference \cite{Mordukhovich2012}, coerciveness is sufficient to guarantee existence. Because $D(\V{x}) \le D_{\epsilon}(\V{x})$, the perturbed criterion $D_{\epsilon}(\V{x})$ is coercive whenever the original criterion $D(\V{x})$ is coercive. Henceforth, we will assume that  $S$ or at least one of the $C_i$ is bounded.

A strictly convex function possesses at most one minimum point on a convex set. The function $|x|$ shows that this sufficient condition for uniqueness is hardly necessary. In the Fermat-Weber problem, where the closed convex sets $C_i=\{\V{x}_i\}$ are singletons, the function $D(\V{x})$ is strictly convex if and only if the points $\V{x}_i$ are non-collinear. To generalize this result, we require the sets $C_i$ to be non-collinear. Geometrically this says that it is impossible to draw a straight line that passes through all of the $C_i$. Non-collinearity can only be achieved when $k>2$ and $\cap_{i=1}^k C_i = \emptyset$. We also require the $C_i$ to be strictly convex.  A set $C$ is said to be strictly convex if the interior of the line segment $[\V{x},\V{y}]$ connecting two different points $\V{x}$ and $\V{y}$ of $C$ lies in the interior of $C$. Put another way, the boundary of $C$ can contain no line segments. A  singleton or a closed ball is strictly convex, but a closed box is not.

\begin{proposition}
If the closed convex sets $C_1,\ldots,C_k$ are strictly convex but not collinear, then $D(\V{x})$ is strictly convex.
\end{proposition}
\begin{proof}
Suppose the contrary is true, and choose $\V{x} \not = \V{y}$ and $\alpha$ strictly between 0 and 1 so that 
\begin{equation}
\label{eq:flat}
D[\alpha \V{x} + (1-\alpha)\V{y}] = \alpha D(\V{x}) + (1-\alpha)D(\V{y}).
\end{equation}
Let $L$ be the line $\{s \V{x}+(1-s) \V{y}: s \in \Real\}$ passing through the points $\V{x}$ and $\V{y}$. Then there exists at least one $C_j$ such that $L \cap C_j = \emptyset$. In particular,  $\V{x}$, $\V{y}$, and $ \alpha \V{x} + (1-\alpha)\V{y}$ all fall outside this $C_j$.  Equality (\ref{eq:flat}) implies that 
\begin{equation*}
\begin{split}
& \; \alpha \lVert \V{x} - P_{C_j}(\V{x}) \rVert 
+  (1-\alpha)\lVert \V{y} - P_{C_j}(\V{y}) \rVert \\
=& \;\; \,\, \lVert \alpha \V{x} + (1-\alpha)\V{y} - P_{C_j}[\alpha \V{x}+(1-\alpha) \V{y} ] \rVert \\
\leq& \;\;\,\, \lVert \alpha \V{x} + (1-\alpha)\V{y} - \alpha P_{C_j}(\V{x}) - (1-\alpha) P_{C_j}(\V{y}) \rVert \\
\leq& \; \alpha \lVert \V{x} - P_{C_j}(\V{x}) \rVert + (1-\alpha)\lVert \V{y} - P_{C_j}(\V{y}) \rVert.
\end{split}
\end{equation*}
Since the projection of a point onto $C_j$ is unique, these sandwich inequalities entail
\begin{equation*}
P_{C_j}[\alpha \V{x}+(1-\alpha) \V{y} ]  = \alpha P_{C_j}(\V{x}) + (1-\alpha) P_{C_j}(\V{y}).
\end{equation*}
If $P_{C_j}(\V{x}) \ne P_{C_j}(\V{y})$, then the strict convexity of $C_j$ implies the convex combination $\alpha P_{C_j}(\V{x}) + (1-\alpha) P_{C_j}(\V{y})$ is interior to $C_j$.  Hence, this point cannot be the closest point to the external point $\alpha \V{x} + (1-\alpha) \V{y}$.  Therefore, consider the possibility $P_{C_j}(\V{x}) = P_{C_j}(\V{y})=\V{z}$.  Equality can occur in the inequality
\begin{equation*}
\lVert \alpha \V{x} + (1-\alpha)\V{y} - \V{z} \rVert \leq \alpha \lVert  \V{x} -\V{z} \rVert+ (1-\alpha) \lVert \V{y} -\V{z}\rVert
\end{equation*}
only when $\V{x} -\V{z} = t(\V{y} -\V{z})$  for some $t \ne 1$.  This relation shows that 
\begin{equation*}
\V{z}  = \frac{1}{1-t} \V{x} -\frac{t}{1-t}\V{y}
\end{equation*}
belongs to $L \cap C_j$, contradicting our hypothesis.  Thus, $D(\V{x})$ is strictly convex. 
\end{proof}

The next result shows that the function $D_{\epsilon}(\V{x})$ inherits strict convexity from $D(\V{x})$. Therefore, when $D(\V{x})$ is strictly convex,
$D_{\epsilon}(\V{x})$ possesses a unique minimum point.

\begin{proposition}
If $D(\V{x})$ is strictly convex, then $D_{\epsilon}(\V{x})$ is also strictly convex.
\end{proposition}
\begin{proof}
Fix arbitrary $\V{x} \not = \V{y}$ and $\alpha$ strictly between 0 and 1. The strict convexity of $D(\V{x})$ implies that there is at least one $j$ such that
\begin{equation*}
d(\alpha \V{x} + (1-\alpha)\V{y}, C_j) < \alpha d(\V{x}, C_j) + (1-\alpha) d(\V{y}, C_j).
\end{equation*}
The strict inequality
\begin{equation*}
\begin{split}
\sqrt{d(\alpha \V{x} + (1-\alpha)\V{y}, C_j)^2 + \epsilon} &< 
\sqrt{[\alpha d(\V{x}, C_j) + (1-\alpha) d(\V{y}, C_j)]^2 + \epsilon}, \\
&\leq \alpha \sqrt{d(\V{x}, C_j)^2 + \epsilon} + (1-\alpha)\sqrt{d(\V{y}, C_j)^2 + \epsilon},
\end{split}
\end{equation*}
follows because the function $f_{\epsilon}(u)=\sqrt{u^2+\epsilon}$ is a strictly increasing and convex. Summing over $j$ gives 
the desired result.
\end{proof}

We now clarify the relationship between the minima of the $D_{\epsilon}(\V{x})$ and $D(\V{x})$ functions.
\begin{proposition}\label{prop:MM_cluster}
 For a sequence of constants 
$\epsilon_{m}$ tending to 0, let $\V{y}_{m}$ be a corresponding 
sequence minimizing $D_{\epsilon_m}(\V{x})$.  If $\V{y}$ is the unique
minimum point of $D(\V{x})$, then $\V{y}_m$ tends to $\V{y}$.  If $D(\V{x})$ has multiple minima, then every cluster point of the sequence $\V{y}_m$ minimizes $D(\V{x})$.
\end{proposition}
\begin{proof}
To prove the assertion, consider the inequalities
\begin{equation*}
D(\V{y}_m) \le D_{\epsilon_m}(\V{y}_{m}) \le D_{\epsilon_m}(\V{x}) \le D_1(\V{x})
\end{equation*}
for any $\V{x} \in S$ and $\epsilon_m \le 1$. Taking limits along the appropriate subsequences
proves that the cluster points of the sequence $\V{y}_m$ minimize $D(\V{x})$. Convergence to a unique minimum point $\V{y}$ occurs provided the sequence $\V{y}_m$ is bounded. If $S$ is bounded, then $\V{y}_m$ is bounded by definition. On the other hand,
if any $C_j$ is bounded, then $D(\V{x})$ is coercive, and the  inequality $D(\V{y}_m) \le D_1(\V{x})$ forces $\V{y}_m$ to be bounded.
\end{proof}

The convergence theory of MM algorithms hinges on the properties of the algorithm map $\psi(\V{x}) \equiv \arg\min_{\V{y}} g(\V{y} \mid \V{x})$. 
For easy reference, we state a simple version of Meyer's monotone convergence theorem \cite{Meyer1976} instrumental in proving convergence in our setting.

\begin{proposition}\label{prop:MM_limit_points}
  Let $f(\V{x})$ be a continuous function on a domain $S$ and
   $\psi(\V{x})$ be a continuous algorithm map from $S$ into $S$ satisfying
 $f(\psi(\V{x})) < f(\V{x})$ for all $\V{x} \in S$ with $\psi(\V{x}) \neq \V{x}$.
  Suppose for some initial point $\V{x}_{0}$ that the set
  $\mathcal{L}_f(\V{x}_{0}) \equiv \{\V{x} \in S : f(\V{x}) \leq f(\V{x}_{0}) \}$ is compact.
Then
  \begin{inparaenum}[(a)]
  \item \label{part:fixed_points} 
    all cluster points are fixed points of $\psi(\V{x})$, and
  \item \label{part:successive_iterates}
     $\lim_{m \to \infty} \lVert \V{x}_{m+1} - \V{x}_{m} \rVert = 0$.
  \end{inparaenum}
\end{proposition}
Note that Proposition~\ref{prop:MM_limit_points} also ensures the existence of at least one cluster point for the sequence of iterates $\V{x}_{m+1} = \psi(\V{x}_{m})$.
Additionally, the convergence of the MM iterates (\ref{eq:MM_update}) to a stationary point of $f(\V{x})$ follows immediately provided the fixed points of $\psi(\V{x})$ are stationary points of $f(\V{x})$ and $\psi(\V{x})$ possesses only finitely many fixed points.
  
Let us verify the conditions of Proposition~\ref{prop:MM_limit_points} for minimizing $D_\epsilon(\V{x})$. The function $D_\epsilon(\V{x})$ is continuous on its domain $S$, and the set $\mathcal{L}_{D_\epsilon}(\V{x}_0)$ is compact for any initial point $\V{x}_0$ since either $S$ is compact or $D_{\epsilon}(\V{x})$ is coercive. The continuity of the algorithm map follows immediately from the continuity of the projection mapping. Finally, we need to prove that $D_{\epsilon}(\psi(\V{x})) < D_{\epsilon}(\V{x})$ whenever $\V{x} \not = \psi(\V{x})$.
First observe that $\psi(\V{x}) = \V{x}$ if and only if the MM surrogate function satisfies $g_{\epsilon}(\V{x} \mid \V{x}) = \min_\V{y} g_{\epsilon}(\V{y} \mid \V{x})$. Since $g_{\epsilon}(\V{y} \mid \V{x})$ has a unique minimizer, we have the strict inequality $g_{\epsilon}(\psi(\V{x}) \mid \V{x}) < g_{\epsilon}(\V{x} \mid \V{x})$ whenever $\V{x}$ is not a fixed point of $\psi$. This forces a decrease in the objective function $D_{\epsilon}(\V{x})$ and makes the MM algorithm strictly monotone outside the set of stationary points.

We now argue that the fixed points of the algorithm map $\psi(\V{x})$ are stationary points of $D_\epsilon(\V{x})$. We will show, in fact, 
that the two sets of points coincide. To accomplish this,
we need to determine the gradients of $D_\epsilon(\V{x})$ and $g_{\epsilon}(\V{x} \mid \V{y})$.
Recall that $f_{\epsilon}(u)$ is strictly increasing and strictly convex.
As a consequence the functions $f_{\epsilon}(\|\V{x}\|)$ and
$f_{\epsilon}[d(\V{x},C_j)]$ are convex.  Even more remarkable is the fact
that both functions are continuously differentiable. When $\V{x} \ne {\bf 0}$, the function $\lVert \V{x} \rVert$ is differentiable. Likewise, when
$\V{x} \not \in C_j$, the function $d(\V{x}, C_j)$ is differentiable. Therefore, the chain rule implies
\begin{eqnarray}
\label{eq:grads_D}
\nabla f_{\epsilon}(\|\V{x}\|) & = & \frac{\|\V{x}\| }{ \sqrt{\|\V{x}\|^2+\epsilon}}\frac{\V{x} }{ \|\V{x}\|}
\:\;\, = \:\;\, \frac{\V{x} }{\sqrt{\|\V{x}\|^2+\epsilon} }\\
\label{eq:grads_g}
\nabla f_{\epsilon}[d(\V{x},C_j)] & = & \frac{d(\V{x},C_j) }{ \sqrt{d(\V{x},C_j)^2+\epsilon}}
\frac{\V{x} - P_{C_j}(\V{x}) }{ d(\V{x},C_j)} \:\;\, = \:\;\, \frac{\V{x} - P_{C_j}(\V{x}) }{ \sqrt{d(\V{x},C_j)^2+\epsilon}},
\end{eqnarray}
respectively. 

 By continuity one expects the gradients to be defined for $\V{x} = {\bf 0}$
and $\V{x} \in C_j$ by the corresponding limit of $\bf 0$.  In the
former case the expansion
\begin{eqnarray*}
\sqrt{\|\V{x}\|^2+\epsilon}-\sqrt{\epsilon} & = & \sqrt{\epsilon} \sqrt{1+\frac{\|\V{x}\|^2}{\epsilon}} - \sqrt{\epsilon}
\:\;\, = \:\;\, \frac{1 }{ 2} \frac{\|\V{x}\|^2}{\sqrt{\epsilon}}+ \sqrt{\epsilon} o\left(\frac{\|\V{x}\|^2}{\epsilon} \right) .
\end{eqnarray*}
shows that $\nabla f_{\epsilon}(\|{\bf 0}\|) = {\bf 0}$.
In the latter case the expansion
\begin{eqnarray*}
\sqrt{d(\V{y},C_j)^2+\epsilon}-\sqrt{\epsilon} & = & \frac{1 }{ 2} \frac{d(\V{y},C_j)^2 }{ \sqrt{\epsilon}}+ 
\sqrt{\epsilon} o\left[\frac{d(\V{y},C_j)^2 }{ \epsilon} \right] 
\end{eqnarray*}
and the bound $d(\V{y},C_j) = |d(\V{y},C_j)-d(\V{x},C_j)| \le \|\V{y}-\V{x}\|$ for $\V{x} \in C_j$
likewise show that $\nabla f_{\epsilon}[d(\V{x},C_j)]= {\bf 0}$. Consequently, equations 
(\ref{eq:grads_D}) and (\ref{eq:grads_g}) hold for all $\V{x} \in \Real^d$. It follows that both $D_{\epsilon}(\V{x})$ and $g_{\epsilon}(\V{x} \mid \V{y})$ are differentiable on $\Real^d$, with gradients
\begin{equation*}
\nabla D_{\epsilon}(\V{x}) \:\;\, = \:\;\, \sum_{j=1}^k 
\gamma_j \frac{\V{x} - P_{C_j}(\V{x}) }{ \sqrt{d(\V{x},C_j)^2+\epsilon}},
\end{equation*}
and
\begin{equation} 
\label{eq:fixed}
\nabla g_{\epsilon}(\V{x} \mid \V{y}) \:\;\, = \:\;\, \sum_{j=1}^k 
\gamma_j \frac{\V{x} - P_{C_j}(\V{y}) }{ \sqrt{d(\V{y},C_j)^2+\epsilon}},
\end{equation}
respectively. Note that $\V{y} \in S$ minimizes $D_{\epsilon}(\V{x})$ over $S$ if and only if
\begin{equation*}
\sum_{j=1}^k 
\gamma_j \frac{\langle \V{y} - P_{C_j}(\V{y}), \V{x} - \V{y} \rangle }{ \sqrt{d(\V{y},C_j)^2+\epsilon}}  \geq 0,
\end{equation*}
for all $\V{x} \in S$. This inequality, however, is equivalent to the inequality $\langle \nabla g_{\epsilon}(\V{y} \mid \V{y}), \V{x} - \V{y}\rangle \geq 0$,
for all $\V{x} \in S$, which in turn holds if and only if $\V{y}$ is a fixed point of $\psi(\V{x})$. If $D(\V{x})$ is strictly convex, then $D_{\epsilon}(\V{x})$ has a unique minimum point, and $\psi(\V{x})$ has exactly one fixed point.

Thus, Proposition~\ref{prop:MM_cluster} and Proposition~\ref{prop:MM_limit_points} together tell us that 
 $\V{y}$ is a solution to (\ref{eq:ourHeron})
if there is a sequence of $\epsilon_m$ tending to zero and a sequence of points $\V{y}_m$ tending to $\V{y}$ that satisfy
\begin{equation}
\left \langle -\sum_{j=1}^k  \gamma_j \frac{ \V{y}_m - P_{C_j}(\V{y}_m) }{ \sqrt{d(\V{y}_m,C_j)^2+\epsilon_m}} , \V{x} - \V{y}_m \right \rangle \leq 0,
\end{equation}
for all $\V{x} \in S$. The above sufficient condition becomes necessary as well if $D(\V{x})$ is strictly convex. As a sanity check, when the sets $S \cap C_j$ are all empty and the weights $\gamma_j$ are identical, we recover the characterization of the optimal points given in Theorem 3.2 of reference \cite{Mordukhovich2012}, albeit under the more restrictive assumption of strict convexity.

\section{Conclusion}

There is admittedly an art to applying the MM principle. The majorization presented here is specific to Euclidean distances, and changing the underlying norm would require radical revision. Nonetheless, when the MM principle applies, the corresponding MM algorithm can be effective, simple to code, and intuitively appealing. Here the principle lit the way to an efficient numerical algorithm for solving the Euclidean version of the generalized Heron problem using only elementary principles of smooth convex analysis.  We also suggested a simple yet accurate approximation of the problem  that removes the singularities of the MM algorithm and Weiszfeld's earlier algorithm. Similar advantages accrue across a broad spectrum of optimization problems.  The ability of MM algorithms to handle high-dimensional problems in imaging, genomics, statistics, and a host of other fields testifies to the potency of a simple idea consistently invoked. Mathematical scientists are well advised to be on the lookout for new applications.

\section*{Acknowledgments.}
This research was supported by the United States Public Health Service
grants GM53275 and HG006139.

\bibliographystyle{monthly}
\bibliography{Heron}

\vspace*{.4in}
\noindent {\em Department of Human Genetics, University of California, Los Angeles, CA 90095 \newline
ecchi@ucla.edu}

\vspace*{.2in}
\noindent {\em Departments of Human Genetics, Biomathematics, and Statistics, University of California, Los Angeles, CA 90095 \newline
klange@ucla.edu}

\end{document}